\documentclass[12pt]{article}

\usepackage{geometry}
\usepackage{fouriernc}
\usepackage{amsmath}
\usepackage[svgnames]{xcolor}
\usepackage{hyperref}
\hypersetup{
    colorlinks=true,
    citecolor=Green,
    urlcolor=Crimson
    }
\usepackage{orcidlink}

\title{Irrational Acceleration\\of a Continued Fraction of $\pi$}
\author{Tomasz Stachowiak\orcidlink{0000-0001-9851-9131}\\
{\small\texttt{tomasz@monodromy.group} }}

\date{}

\begin{document}


\maketitle

\begin{abstract}
\noindent An application of (iterated) Bauer-Muir acceleration can give an Ap\'ery-like continued fraction for $\pi$ with irrational coefficients, and much faster convergence. It can be considered a generalized continued fraction with the same matrix representation as the standard case, but the dimension increased due to irrationality. The construction is also given for $\ln(2)$ and $\sqrt[3]{2}$.
\end{abstract}

\section{}

In a recent paper \cite{Cohen:2024}, Henri Cohen described a method of accelerating the convergence of continued fractions, and how its iterative application leads to Ap\'ery-like fractions.

A limitation indicated in that article is that for fractions where irrational number appear in the rate of convergence, one cannot find exact (rational) expressions required to construct the recurrent hierarchy of new continued fractions. Specifically, Cohen writes that when the convergence is of the form $\mathcal{O}(B^n)$, where $B$ is an irrational number, the Bauer-Muir iteration is not applicable.

However, the procedure still works in principle, and, if one allows certain degree of irrationality, leads to practical results. An example of particular interest is the acceleration of the continued fraction of $\pi$, given here. It leads to a generalized fraction, still easily computable despite the presence of an algebraic number.
Such approximations have a better optimal bound on their error, given by the Davenport-Schmidt theorem \cite{Davenport}.

Although there are series expressions for $\pi$ involving roots of integers in their terms, like the Ramanujan or Chudnovsky series, there seems to be no expansion based on $\sqrt{2}$, of the form obtained here.

\section{}

The starting point is the continued fraction of $\pi$ obtained from the arc tangent:
\begin{equation}
    \pi = \cfrac{4}{
        1+\cfrac{1^2}{
            3 + \cfrac{2^2}{
                5 + \cfrac{3^2}{
                    7 + \raisebox{-.33\height}{$\ddots$}\vphantom{\ddots}
                }
            }
        }
    }
\end{equation}
Using the notation of \cite{Cohen:2024}, the above can be rewritten as
\begin{equation}
    \frac{4}{\pi} = a(0) + \cfrac{b(0)}{
        a(1) + \cfrac{b(1)}{
            a(2) + \raisebox{-.33\height}{$\ddots$}\vphantom{\ddots}
        }
    },
\end{equation}
where $a(n) = 2n+1$ and $b(n) = (n+1)^2$. It is straightforward to check the rate of convergence using Theorem 7.3.7 (p. 306) of \cite{Belabas}, which gives
\begin{equation}
    \mathcal{O}\left( (-1)^n \left(1+\sqrt{2}\right)^{-2n} \right)
    \approx \mathcal{O}\left( (-5.83)^n \right).
\end{equation}
The irrationality means that the fraction cannot be accelerated, and the specific reason is that the resulting expansion would have non-integer elements -- it could not be considered a standard continued fraction.

However, in view of both the computations required to complete the procedure and the subsequent computation of the fraction, this does not create major obstacles. The construction can be carried out in the quadratic extension
$\mathbb{Q}\left(\sqrt{2}+1\right)$, leading to a generalized continued fraction, as follows.

The Bauer-Muir acceleration step consists in finding a function $r(n)$ which satisfies the non-linear recurrence
\begin{equation}
    0 = r(n) \big( a(n+1) + r(n+1)\big) - b(n).
\end{equation}
Because $r(n)$ is in fact the $n$-th tail of the continued fraction, there is little hope of finding it explicitly as a simple function of $n$. We can, nevertheless, approximate it by trying to solve a modified recurrence
\begin{equation}
    d(n) = r(n) \big( a(n+1) + r(n+1)\big) - b(n),
\end{equation}
where $d(n)$ should be minimal. This means, in practice, that it should be a small constant, or a polynomial of as low degree as possible.

The problem with $\pi$ is that the convergence rate involves irrational numbers, and this translates directly into the recurrence. The optimal solution is
\begin{equation}
\begin{aligned}
    r(n) &= \left(\sqrt{2}-1\right)\left(n+\frac12\right),\\
    d(n) &= -\frac14,
\end{aligned}
\end{equation}
from which the accelerated partial numerators $a(n,1)$ and denominators $b(n,1)$ can be calculated according to
\begin{equation}
\begin{aligned}
    a(n,1) &= a(n+1) + r(n+1) &\text{for} && n < 1,\\
    a(n,1) &= a(n+1) + r(n+1) -r(n-1)\frac{d(n)}{d(n-1)}, &\text{for} && n\geq 1,\\
    b(-1,1) &= -d(0),\\
    b(n,1) &= b(n)\frac{d(n+1)}{d(n)} &\text{for} && n \geq 0.
\end{aligned}
\end{equation}
These are the same as in \cite{Cohen:2024}, except for reindexing, which was introduced for convenience: the generic recurrence starts at 0 with $b(0,1)$.
In the present case, the explicit expressions are
\begin{equation}
\begin{aligned}
    a(-1,1) &= \frac12\left(\sqrt{2}+1\right),\\
    a(0,1) &= \frac32\left(\sqrt{2}+1\right),\\
    a(n,1) &= 2n+1+2\sqrt{2} &\text{for}\quad n\geq 1,\\
    b(-1,1) &= \frac14,\\
    b(n,1) &= (n+1)^2 &\text{for} \quad n\geq 0.
\end{aligned}
\end{equation}

At this stage, the whole fraction is of the form
\begin{equation}
    \frac{4}{\pi} = a(-1,1) + \cfrac{b(-1,1)}{
        a(0,1)+\cfrac{b(0,1)}{
            a(1,1)+\raisebox{-.33\height}{\(\ddots\)}\vphantom{\ddots}
        }
    },
\end{equation}
and the whole step can be applied again, yielding $a(n,2)$ from $a(n,1)$ just as
$a(n,1)$ was obtained from $a(n)$, and similarly for $b$. The Ap\'ery acceleration combines infinitely many such steps, giving a fraction whose terms are obtained by incrementing both the depth of the fraction and the Bauer-Muir index -- these are the two indices $n$ and $l$ in $a(n,l)$ and $b(n,l)$.

More precisely, we have to find a system of generic functions $r(n,l)$, $d(n,l)$, $a(n,l)$, and $b(n,l)$, which satisfy the mutual recurrence equations
\begin{equation}
\begin{aligned}
    d(n,l) &= r(n,l) R(n+1,l) - b(n,l),\\
    a(n,l+1) &= R(n+1,l) - r(n-1,l)\frac{d(n,l)}{d(n-1,l)},\\
    b(n,l+1) &= b(n,l)\frac{d(n+1,l)}{d(n,l)},\\
    R(n,l) &= a(n,l) + r(n,l),\\
\end{aligned}
\end{equation}
where the last function, $R$, is merely a convenient shorthand. The word `generic' reflects the fact, that the initial terms might require a separate calculation. This point will be addressed below.

By inspection of several explicit Bauer-Muir steps ($l$ up to about 4) it is possible to guess and verify the the sought-for functions. It luckily turns out that no other irrational numbers appear in the subsequent stages, and the general procedure can be continued. It will thus be convenient to denote the special algebraic number $\sqrt{2}+1$ as $\omega$ for brevity. Its minimal polynomial is
\begin{equation}
    \mu(\omega) = \omega^2 -2\omega -1,
\end{equation}
and it can be employed to reduce all rational expression involving $\omega$ through
\begin{equation}
    \omega^2 = 2\omega+1, \quad \text{and} \quad 
    \frac{1}{v+\omega z} =  \frac{v+2z-\omega z}{v^2+2 v z-z^2}.
\end{equation}
This makes the calculations with symbolic packages like Mathematica much more efficient than using $\sqrt{2}$ explicitly. Keep in mind also, that the choice is not unique: $\sqrt{2}$ would do as well at first, but some of the polynomial expressions below would be much longer. The next section gives further possible simplifications tied to the choice of $\omega$.

The generic functions can now be given as
\begin{equation}
\begin{aligned}
    r(n,l) &= (\omega-2)\left(n-l+\frac12\right)\\
    d(n,l) &= -\frac14(2l+1)^2,\\
    a(n,l) &= 2n + 1 + 2(\omega-1) l,\\
    b(n,l) &= (n+1)^2.
\end{aligned}
\end{equation}
This allows us to construct the partial numerators $p(n,n)$ and denominators $q(n,n)$ that correspond to a diagonal walk in the index space: the $n$-th partial expression after $n$ Bauer-Muir steps. 
Just as in the case of ordinary continued fractions there is a two-term recurrence for both $p(n)$ and $q(n)$, we have two-term recurrences for $p(n,l)$ and $q(n,l)$: one for increasing each index. 
They can be found in Corollary 7.5.3 of \cite{Belabas}:
\begin{equation}
\begin{aligned}
    u(n+1,l+1) &= R(n+2,l) u(n,l+1) - d(n+1,l)u(n,l),\\
    u(n,l+1) &= R(n+1,l) u(n,l) + b(n,l) u(n-1,l),
\end{aligned}
\end{equation}
where $u$ stands for either $p$ or $q$.
This is perhaps clearer when written as the explicit continued fraction (Corollary 7.5.4 of \cite{Belabas})
\begin{equation}
    \frac{p(n,n)}{q(n,n)} = 
    a(0) + \cfrac{b(0)}{
        R(1,0) + \cfrac{-d(1,0)}{
            R(2,0) + \cfrac{b(1,1)}{
                R(2,1) + \cfrac{-d(2,1)}{
                    R(3,1) + \vphantom{\ddots}\raisebox{-.33\height}{\(\,\ddots\,\)}
                    \raisebox{-.4\height}{\( + \cfrac{-d(n,n-1)}{R(n+1,n-1)}\)
                    }
                }
            }
        }
    }.
    \label{diag_conv}
\end{equation}
This is a period 2 fraction: the partial expressions alternate, and it would be convenient to contract them into one numerator/denominator pair. In other words,
to transform the present form
\begin{equation}
    W_0(n) = R(n+1,n) + \cfrac{-d(n+1,n)}{
            R(n+2,n) + \cfrac{b(n+1,n+1)}{W_0(n+1)}
        }
        \label{init_W}
\end{equation}
into the familiar
\begin{equation}
    W(n) = \frac{P(n)}{Q(n)+W(n+1)}.
    \label{simple_W}
\end{equation}
This can be done in two steps. First, by shifting:
\begin{equation}
    W_0(n) = W_1(n) + 4n+3 + \frac{2(n+1)(5n+4)}{\omega(4n+5)},
\end{equation}
and then by rescaling:
\begin{equation}
    W_1(n) = \frac{W_2(n)}{\omega(1+4n)(5+4n)}.
\end{equation}
These can simply be determined from writing out \eqref{init_W} explicitly, and first deleting the free term to obtain a single fraction of the form \eqref{simple_W},
then getting rid of the multiplicative coefficient in front of $W$. The result of the substitution in \eqref{init_W} is
\begin{equation}
    W_2(n) = \frac{-(k+2)^2(2k+1)^2(4k+1)(4k+9)}{
        6(k+2)(2k+3)(4k+7) + 
        (4k+5)(4k+7)(4k+9)\omega + W_2(n+1)}.
        \label{gen_rec_W}
\end{equation}

Finally, the initial terms have to be determined by hand, although thanks to keeping track of all transformation it is not necessary to use the numerical method of \cite{Cohen:2024}. We can rewrite \eqref{diag_conv} with its first terms up to $W_0$, which can then be expressed through $W_2$ thus:
\begin{equation}
    a(0) + \cfrac{b(0)}{W_0(0)} = 
    1 + \cfrac{5\omega}{8+15\omega + W_2(0)}.
\end{equation}
Starting from $W_2(0)$, the generic recurrence \eqref{gen_rec_W} holds, so we arrive at
the complete (generalized) continued fraction
\begin{equation}
    \frac{4}{\pi} = 
    1 + \cfrac[l]{\quad\; 5\omega}{8+15\omega + \cfrac[l]{\quad\; S_0}{
    T_0+\omega\,U_0 + \cfrac[l]{\quad\; S_1}{
        T_1 + \omega\, U_1 + \raisebox{-.33\height}{\(\ddots\)}\vphantom{\ddots}
    }
    }}
    \label{result}
\end{equation}
where
\begin{equation}
\begin{aligned}
S_k &= -(k+2)^2(2k+1)^2(4k+1)(4k+9),\\
T_k &= 6(k+2)(2k+3)(4k+7),\\
U_k &= (4k+5)(4k+7)(4k+9).
\end{aligned}
\label{STU_polys}
\end{equation}

Using Theorem 7.3.7\cite{Belabas} again, the rate of convergence can be checked to be
\begin{equation}
    \mathcal{O}\left( \left(7+4\sqrt{2}+2\sqrt{20+14\sqrt{2}}\right)^{-2n}\right) 
    \approx 639^{-n}
    \label{improved}
\end{equation}
-- a considerable improvement, even though the type of convergence has not changed. The new fraction gives 2.8 decimal places per step, where the original one gave 0.77.

The question remains as to how to deal with the square root in practical computation, and whether a simple recurrence can be given to obtain successive convergents. Thankfully we can proceed like with complex numbers: treat irrational expression $x+\omega y$ as pairs of integers. This will increase the dimension by 2, but everything will stay linear.

Recall, that for a continued fraction with partial numerators $A_n$, and denominators $B_n$, there is just one recurrence
\begin{equation}
    u_n = B_n u_{n-1} + A_n u_{n-2},
    \label{basic}
\end{equation}
where $u_n$ stands for either $p_n$ or $q_n$, and the $n$-th convergent is $p_n/q_n$. If these numbers are decomposed with respect to $\omega$ as
$p_n = x_n+\omega y_n$, the recurrence for $p_n$ becomes
\begin{equation}
\begin{aligned}
    x_n + \omega y_n &= (T_n+\omega\,U_n)(x_{n-1}+\omega y_{n-1}) +
        S_n (x_{n-2}+\omega y_{n-2}),
\end{aligned}
\end{equation}
or
\begin{equation}
\begin{aligned}
x_n &= T_n x_{n-1} + U_n y_{n-1} + S_n x_{n-2},\\
y_n &= T_n y_{n-1} + U_n x_{n-1} +2U_n y_{n-1} + S_n y_{n-2}.
\end{aligned}
\end{equation}
The same holds with the decomposition of $q_n = v_n+\omega z_n$, and by introducing the matrix
\begin{equation}
    \mathbf{Y}_n = \begin{pmatrix}
    y_{n} & x_{n}\\ 
    z_{n} & v_{n}
    \end{pmatrix},
\end{equation}
we can write both recurrences together
\begin{equation}
    \mathbf{Y}_n 
    = \mathbf{Y}_{n-2}S_n + \mathbf{Y}_{n-1}\left(T_n + \Omega\, U_n\right),
    \label{recur2}
\end{equation}
where
\begin{equation}
    \Omega = \begin{pmatrix}
    2 & 1 \\
    1 & 0
    \end{pmatrix}.
\end{equation}
Notice, that $\Omega$ has $\mu$ as it's characteristic polynomial
$\mu(\Omega) = 0$ -- it is the matrix equivalent of our algebraic number $\omega$. We have recovered the exact form of \eqref{basic} without any irrational numbers, at the cost of passing to matrices.

One further simplification is possible: upon defining block matrices
\begin{equation}
    \mathbf{X} _n = \begin{pmatrix}
    \mathbf{Y}_{n-1} & \mathbf{Y}_{n}
    \end{pmatrix}, \quad
    \mathbf{M} _n = \begin{pmatrix}
        0 & S_n \\
        1 & T_n + \Omega\, U_n
    \end{pmatrix},
\end{equation}
the recurrence becomes the first-order system
\begin{equation}
    \mathbf{X}_n = \mathbf{X}_{n-1} \mathbf{M}_{n}.
\end{equation}
In effect, the computation requires keeping track of four integers per convergent instead of two.
The square root can be evaluated at the very end to the desired precision, as the convergents are always rational in $\omega$, but can trivially be made linear through
\begin{equation}
    \frac{p_n}{q_n} = \frac{x_n +\omega\,y_n}{v_n+\omega\,z_n}
    = \frac{x_n v_n + (2x_n-y_n) z_n + \omega\,(v_n y_n - x_n z_n)}
    {v_n^2+2v_n z_n-z_n^2}.
\end{equation}

In implementation, the square root of 2 could itself be approximated by its own continued fraction, or Gosper's algorithm \cite{Gosper} could be used to directly evaluate the whole rational expression.

The whole process could be repeated again, but with new algebraic numbers appearing in \eqref{improved}, another extension would be needed, and this time the minimal polynomial is of the fourth degree. So although convergence would improve, the increase in dimensionality (number of integer components) would slow down the actual computation.

\section{}
There are several continued fractions with simple irrational convergence rates mentioned in \cite{Cohen:2024}, and for some of them an equally simple result exists. The procedure is exactly as above, so only the final formulae are given here.

Just as $\pi = 4\arctan(1)$, we can write $\log(2) = 2\mathrm{artanh}(\frac13)$, and use the same general fraction, which has the rate of convergence $\omega^{-4n}$,
giving about 1.53 decimal places per step. The generalized continued fraction, written in the notation of Gauss, is of the same form as \eqref{result}:
\begin{equation}
    \ln(2) = \cfrac{6\eta}{
        -2 + 9\eta - 
        \underset{k=0}{\overset{\infty}{\text{{\large K}}}}\,\cfrac{S_k}{T_k+\eta\,U_k}
    },
\end{equation}
with $\eta = \omega^2 = 3+2\sqrt{2}$, and 
\begin{equation}
\begin{aligned}
S_k &= -(1+k)^2(1+2k)^2(-1+4k)(7+4k),\\
T_k &= -2(5+4k)(8+15k+6k^2),\\
U_k &= 3(3+4k)(5+4k)(7+4k).
\end{aligned}
\end{equation}
We notice that $T_k$ is no longer fully factorized, and that its form depends on the choice of the algebraic number $\eta$ -- the partial denominator could just as well be written as $T_k+ cU_k + (\eta-c)U_k$ for some rational $c$.

In the present case, the first few simple rational $c$ can be checked by hand, and it turns out, that $c=\frac14$ simplifies the whole fraction considerably.
We have $4T_k+U_k = -(5+4k)$, and so the fraction becomes
\begin{equation}
    \ln(2) = \cfrac{6+6\eta_1}{
        1+9\eta_1 - 
        \underset{k=0}{\overset{\infty}{\text{{\large K}}}}
        \,\cfrac{S'_k}{T'_k+\eta'\,U'_k}
    },
\end{equation}
with $\eta' = 4\eta-1$, and
\begin{equation}
\begin{aligned}
S'_k &= -16(1+k)^2(1+2k)^2(-1+4k)(7+4k),\\
T'_k &= -(5+4k),\\
U'_k &= 3(3+4k)(5+4k)(7+4k).
\end{aligned}
\end{equation}
This saves two multiplications at each step. The final rate of convergence is
\begin{equation}
\left(33+24\sqrt{2}+4\sqrt{140+99\sqrt{2}}\right)^{-2n}
\approx 10^{-4.25n}.
\label{huh}
\end{equation}

The next example is  $\sqrt[3]{2}$, for which several acceleration steps in \cite{Cohen:2024} give the fraction after two full accelerations. There appears to be a minor sign error in \cite{Cohen:2024}, but the final expansion is
\begin{equation}
    \sqrt[3]{2} = -\cfrac12 - \cfrac{4}{
        5 + \underset{k=0}{\overset{\infty}{\text{\large K}}}
        \;\cfrac{16-9k^2}{9(2k+1)}
    },
\end{equation}
and it too converges like $\omega^{-4n}$. The additional irrational acceleration yields
\begin{equation}
    \sqrt[3]{2} = \cfrac12 +\cfrac{108\eta}{
        -288 + 108\eta + \underset{k=0}{\overset{\infty}{\text{\large K}}}
        \,\cfrac{S_k}{T_k+\eta\,U_k}
    },
\end{equation}
with
\begin{equation}
\begin{aligned}
S_k &= 9(16-9k^2)(4k-3)(4k+5)(6k-5)(6k+11),\\
T_k &= -6(4k+3)(54k^2+81k-14),\\
U_k &= 81(4k+1)(4k+3)(4k+5).
\end{aligned}
\end{equation}
Unfortunately, no simpler $T_k$ could be found with a simple rational modification of $\eta$ as before. The rate of convergence is improved, exactly as above, to \eqref{huh}.

Let us close with the final, small improvement for $\pi$. If other algebraic numbers are checked, as for $\ln(2)$, the polynomial $T_k$ can be made linear, again saving two multiplications. Namely, for $\omega' = 4\omega + 3 = 4\sqrt{2}+7$,
and new polynomials
\begin{equation}
\begin{aligned}
S'_k &= -16(2+k)^2(1+2k)^2(1+4k)(9+4k),\\
T'_k &= 9(7+4k),\\
U'_k &= (5+4k)(7+4k)(9+4k),
\end{aligned}
\end{equation}
we have
\begin{equation}
    \pi = \cfrac{4}{
        1+\cfrac{15-5\omega'}{
            13-15\omega' - \underset{k=0}{\overset{\infty}{\text{\large K}}}
        \,\cfrac{S'_k}{T'_k+\omega'\,U'_k}
        }
    }.
\end{equation}

\end{document}